\documentclass[reqno]{amsart}
\usepackage{amssymb}
\usepackage{amsfonts}

\newtheorem{theorem}{Theorem}
\theoremstyle{plain}

\newtheorem{remark}{Remark}

\numberwithin{equation}{section}
\input{tcilatex}

\begin{document}
\title[Easier Proofs For Kato Estimates]{Elementary Proofs for Kato
Smoothing Estimates of Schr\"{o}dinger-like Dispersive Equations}
\author{Xuwen Chen}
\address{Department of Mathematics, University of Maryland, College Park, MD
20742}
\email{chenxuwen@math.umd.edu}
\thanks{}
\date{07/12/2010}

\begin{abstract}
In this expository note, we consider the dispersive equation:%
\begin{equation*}
i\phi _{t}=(-\triangle )^{\frac{\beta }{2}}\phi \text{ }in\text{ }\mathbb{R}%
^{n+1},\text{ }\phi (x,0)=f(x)\in L^{2}(\mathbb{R}^{n}).
\end{equation*}%
We prove some extensions and refinements of classical Kato type estimates
with elementary techniques.
\end{abstract}

\maketitle

In this short note, we give easier and unified proofs for certain smoothing
estimates of the dispersive equation:%
\begin{equation}
i\phi _{t}=(-\triangle )^{\frac{\beta }{2}}\phi \text{ }in\text{ }\mathbb{R}%
^{n+1},\text{ }\phi (x,0)=\phi _{0}(x)\in L^{2}(\mathbb{R}^{n}).
\label{equation}
\end{equation}%
Theorems \ref{FullKatoEstimateForFreeSchrodinger} and \ref{1d odd estimate
for Schrodinger} extend the classical Kato estimate:%
\begin{equation}
\dint\limits_{-\infty }^{\infty }\dint\limits_{\mathbb{R}^{n}}\frac{%
\left\vert \left\vert \nabla \right\vert ^{\alpha }\phi (x,t)\right\vert ^{2}%
}{\left\vert x\right\vert ^{2-2\alpha }}dxdt\leqslant C\left\Vert \phi
(\cdot ,0)\right\Vert _{2}^{2},\text{ for }\alpha \in \lbrack 0,\frac{1}{2})%
\text{ and }n\geqslant 3  \label{KlainermanKato}
\end{equation}%
in Kato and Yajima \cite{Kato}, and Ben-Artzi and Klainerman \cite%
{Klainerman} for the free Schr\"{o}dinger equation $(\beta =2$ in equation %
\ref{equation}) and show that estimate \ref{KlainermanKato} is in fact an
identity whenever the initial data is radial. In particular, this also
implies that when $n=3$ and $\alpha =0,$ the best constant in estimate \ref%
{KlainermanKato} is attained for every $L^{2}(\mathbb{R}^{3})$ radial data
(via Simon \cite{Simon}). As pointed out in Vilela \cite{Vilela}, the free
Schr\"{o}dinger endpoint Strichartz estimate for radial data in the case
when $n\geqslant 3$ follows from estimate \ref{KlainermanKato}. Moreover,
the proof of theorem \ref{1d odd estimate for Schrodinger} in fact gives
theorem \ref{1/2derivative for 1d} which is stated below.

\subsection{Statement of the theorems}

\begin{theorem}
\label{FullKatoEstimateForFreeSchrodinger}Let $\phi $ be the solution to
equation \ref{equation}, then for $1<\beta -2\alpha <n$, we have 
\begin{equation*}
\int_{-\infty }^{\infty }\int_{\mathbb{R}^{n}}\frac{\left\vert \left\vert
\nabla \right\vert ^{\alpha }\phi (x,t)\right\vert ^{2}}{\left\vert
x\right\vert ^{\beta -2\alpha }}dxdt\leqslant C_{n,\alpha ,\beta }\left\Vert
\phi _{0}\right\Vert _{2}^{2},
\end{equation*}%
Moreover, if $\phi _{0}$ is spherically symmetric, then equality holds i.e. 
\begin{equation*}
\int_{-\infty }^{\infty }\int_{\mathbb{R}^{n}}\frac{\left\vert \left\vert
\nabla \right\vert ^{\alpha }\phi (x,t)\right\vert ^{2}}{\left\vert
x\right\vert ^{\beta -2\alpha }}dxdt=C_{n,\alpha ,\beta }\left\Vert \phi
_{0}\right\Vert _{2}^{2}.
\end{equation*}
\end{theorem}

\begin{remark}
As mentioned before, the above estimate when $\beta =2$, was proved by Kato
and Yajima \cite{Kato} in 1989, Ben-Artzi and Klainerman \cite{Klainerman}
in 1992. The case $\beta =2$, $\alpha =0$ was also mentioned by Herbst \cite%
{Herbst} and Simon \cite{Simon} in 1991. Vilela reproved estimate \ref%
{KlainermanKato} to give the endpoint Strichartz estimate for radial data in
the case when $n\geqslant 3$ in \cite{Vilela} in 2001 . However, they did
not show the equality for radial data. In addition, we will avoid the use of
trace lemmas.
\end{remark}

When $n=1,$ we have the same theorem back, but we have to assume odd initial
data.

\begin{theorem}
\label{1d odd estimate for Schrodinger}Let $\phi $ be the solution to
equation \ref{equation}$\ $in $\mathbb{R}^{1+1}$ with odd initial data, i.e.%
\begin{equation*}
\phi _{0}(-x)=-\phi _{0}(x),
\end{equation*}%
then for $1<\beta -2\alpha \leqslant 2$, we have the identity 
\begin{equation*}
\int_{-\infty }^{\infty }\int_{-\infty }^{\infty }\frac{\left\vert
\left\vert \nabla \right\vert ^{\alpha }\phi (x,t)\right\vert ^{2}}{%
\left\vert x\right\vert ^{\beta -2\alpha }}dxdt=\frac{2^{\beta -2\alpha
}\Gamma (2-\beta +2\alpha )\sin \left( \frac{2-\beta +2\alpha }{2}\pi
\right) }{\beta (\beta -1-2\alpha )}\left\Vert \phi _{0}\right\Vert _{2}^{2}.
\end{equation*}%
In particular, when $\alpha =0,$ $\beta =2,$ we have%
\begin{equation*}
\int_{-\infty }^{\infty }\int_{-\infty }^{\infty }\frac{\left\vert \phi
(x,t)\right\vert ^{2}}{\left\vert x\right\vert ^{2}}dxdt=\pi \left\Vert \phi
_{0}\right\Vert _{2}^{2}.
\end{equation*}%
Or equivalently, say $\psi (\left\vert x\right\vert ,t)$ solves equation \ref%
{equation} when $\beta =2$ in $\mathbb{R}^{3+1}$ as a 3d radial function,
then we have the identity 
\begin{equation*}
\int_{-\infty }^{\infty }\int_{\mathbb{R}^{3}}\frac{\left\vert \psi
(\left\vert x\right\vert ,t)\right\vert ^{2}}{\left\vert x\right\vert ^{2}}%
dxdt=\pi \left\Vert \psi (\left\vert \cdot \right\vert ,0)\right\Vert
_{L^{2}(\mathbb{R}^{3})}^{2}.
\end{equation*}
\end{theorem}

\begin{remark}
Simon showed that the best constant in the classical Kato estimate \ref%
{KlainermanKato} is $\frac{\pi }{n-2}$ when $\alpha =0$ in \cite{Simon}, but
he did not give an explicit $\phi _{0}$ to reach that bound.
\end{remark}

\begin{remark}
It is true that if 
\begin{equation}
iu_{t}=-\triangle u+\left\vert x\right\vert ^{2}u\text{ in }\mathbb{R}^{n+1}%
\text{,}  \label{quadratic Schrodinger}
\end{equation}%
then%
\begin{equation*}
\dint\limits_{0}^{2\pi }\dint\limits_{\mathbb{R}^{n}}\frac{\left\vert
u(x,t)\right\vert ^{2}}{\left\vert x\right\vert ^{2}}dxdt\leqslant
C\left\Vert u(\cdot ,0)\right\Vert _{2}^{2}
\end{equation*}%
when $n\geqslant 3$. Also there is a theorem similar to theorem \ref{1d odd
estimate for Schrodinger} for equation \ref{quadratic Schrodinger} in $%
\mathbb{R}^{1+1}.$ However, the proof is quite different from what we are
dealing with here. See Chen \cite{Chen DIE}.
\end{remark}

\begin{remark}
Vega and Visciglia also proved a family of identities involving the local
smoothing effect for the Schr\"{o}dinger equation. See Vega and Visciglia 
\cite{VegaAndVisciglia}.
\end{remark}

For $\alpha =\frac{\beta -1}{2}$ and $n=1$, the proof of theorem \ref{1d odd
estimate for Schrodinger} in fact reproduces the following result which was
part of theorem 4.1 in Kenig, Ponce and Vega \cite{KenigPonceVega}.

\begin{theorem}
\label{1/2derivative for 1d}Without assuming odd initial data, if $\phi
(x,t) $ solves equation \ref{equation} in $\mathbb{R}^{1+1}$, then we have%
\begin{equation*}
\sup_{x\in \mathbb{R}}\int_{-\infty }^{\infty }\left\vert \left\vert \nabla
\right\vert ^{\frac{\beta -1}{2}}\phi (x,t)\right\vert ^{2}dt\leqslant
C\left\Vert \phi _{0}\right\Vert _{L^{2}(\mathbb{R})}^{2},\text{ }\beta >-1.
\end{equation*}
\end{theorem}

\begin{remark}
The above estimate answers exercise 2.56 in Tao \cite{Tao}.
\end{remark}

\subsection{Proof of theorem \protect\ref{FullKatoEstimateForFreeSchrodinger}%
}

It is well known that%
\begin{equation*}
\left\vert \nabla \right\vert ^{\alpha }\phi (x,t)=\frac{1}{\left( 2\pi
\right) ^{n}}\int_{\mathbb{R}^{n}}\left\vert \xi \right\vert ^{\alpha
}e^{ix\cdot \xi }e^{-i\left\vert \xi \right\vert ^{\beta }t}\hat{\phi}%
_{0}(\xi )d\xi ,
\end{equation*}%
if we choose%
\begin{equation*}
\hat{f}(\xi )=\int_{\mathbb{R}^{n}}e^{-ix\cdot \xi }f(x)dx,
\end{equation*}%
which gives%
\begin{equation*}
\left\Vert \hat{f}\right\Vert _{2}^{2}=\left( 2\pi \right) ^{n}\left\Vert
f\right\Vert _{2}^{2}\text{ and }\int_{\mathbb{R}^{n}}e^{-ix\cdot \xi
}dx=\left( 2\pi \right) ^{n}\delta (\xi ).
\end{equation*}

Hence we have 
\begin{eqnarray}
&&\int_{-\infty }^{\infty }\left\vert \left\vert \nabla \right\vert ^{\alpha
}\phi (x,t)\right\vert ^{2}dt  \label{integrate in time} \\
&=&\frac{1}{(2\pi )^{2n}}\int_{-\infty }^{\infty }dt\int_{\mathbb{R}%
^{n}}d\xi _{1}\int_{\mathbb{R}^{n}}d\xi _{2}\left( e^{ix\cdot \xi
_{1}}e^{-i\left\vert \xi _{1}\right\vert ^{\beta }t}e^{-ix\cdot \xi
_{2}}e^{i\left\vert \xi _{2}\right\vert ^{\beta }t}\left\vert \xi
_{1}\right\vert ^{\alpha }\left\vert \xi _{2}\right\vert ^{\alpha }\hat{\phi}%
_{0}(\xi _{1})\overline{\hat{\phi}_{0}(\xi _{2})}\right)  \notag \\
&=&\frac{1}{(2\pi )^{2n}}\int_{\mathbb{S}^{n-1}}dS_{\omega _{1}}\int_{%
\mathbb{S}^{n-1}}dS_{\omega _{2}}\int_{0}^{\infty }r_{1}^{n-1}dr_{1}  \notag
\\
&&\lim_{\epsilon \rightarrow 0}\int_{0}^{\infty
}r_{2}^{n-1}dr_{2}\int_{-\infty }^{\infty }dt\left( e^{ix\cdot (r_{1}\omega
_{1}-r_{2}\omega _{2})}\eta (\varepsilon t)e^{-i(r_{1}^{\beta }-r_{2}^{\beta
})t}\left\vert r_{1}\right\vert ^{\alpha }\left\vert r_{2}\right\vert
^{\alpha }\hat{\phi}_{0}(r_{1}\omega _{1})\overline{\hat{\phi}%
_{0}(r_{2}\omega _{2})}\right)  \notag \\
&=&\frac{1}{(2\pi )^{2n}}\int_{\mathbb{S}^{n-1}}dS_{\omega _{1}}\int_{%
\mathbb{S}^{n-1}}dS_{\omega _{2}}\int_{0}^{\infty }r_{1}^{n-1}dr_{1}  \notag
\\
&&\lim_{\epsilon \rightarrow 0}\int_{0}^{\infty }r_{2}^{n-1}dr_{2}\left(
e^{ix\cdot (r_{1}\omega _{1}-r_{2}\omega _{2})}\hat{\eta}_{\varepsilon
}(r_{1}^{\beta }-r_{2}^{\beta })\left\vert r_{1}\right\vert ^{\alpha
}\left\vert r_{2}\right\vert ^{\alpha }\hat{\phi}_{0}(r_{1}\omega _{1})%
\overline{\hat{\phi}_{0}(r_{2}\omega _{2})}\right)  \notag \\
&=&\frac{1}{(2\pi )^{2n}}\int_{\mathbb{S}^{n-1}}dS_{\omega _{1}}\int_{%
\mathbb{S}^{n-1}}dS_{\omega _{2}}\int_{0}^{\infty }r_{1}^{n-1}dr_{1}  \notag
\\
&&\lim_{\epsilon \rightarrow 0}\int_{0}^{\infty }v^{\frac{n-1}{\beta }}v^{%
\frac{1}{\beta }-1}\frac{dv}{\beta }\left( e^{ix\cdot (r_{1}\omega _{1}-v^{%
\frac{1}{\beta }}\omega _{2})}\hat{\eta}_{\varepsilon }(r_{1}^{\beta
}-v)\left\vert r_{1}\right\vert ^{\alpha }v^{\frac{a}{\beta }}\hat{\phi}%
_{0}(r_{1}\omega _{1})\overline{\hat{\phi}_{0}(v^{\frac{1}{\beta }}\omega
_{2})}\right)  \notag \\
&=&\frac{1}{\beta }\frac{1}{(2\pi )^{n}}\int_{\mathbb{S}^{n-1}}dS_{\omega
_{1}}\int_{\mathbb{S}^{n-1}}dS_{\omega _{2}}\int_{0}^{\infty }\left(
r_{1}^{n-\beta +2\alpha }e^{ix\cdot (r_{1}\omega _{1}-r_{1}\omega _{2})}\hat{%
\phi}_{0}(r_{1}\omega _{1})\overline{\hat{\phi}_{0}(r_{1}\omega _{2})}%
\right) r_{1}^{n-1}dr_{1}  \notag
\end{eqnarray}%
where $\eta $ is a suitable bump function i.e. $\hat{\eta}_{\varepsilon
}(\xi )=\frac{1}{\varepsilon }\hat{\eta}\left( \frac{\xi }{\varepsilon }%
\right) $ is an approximation to $(2\pi )^{n}\delta (\xi )$. This
approximation of identity is used in order to avoid $\delta (r_{1}^{\beta
}-r_{2}^{\beta })$ in some dimensions.

Whence%
\begin{eqnarray*}
&&\int_{-\infty }^{\infty }\int_{\mathbb{R}^{n}}\frac{\left\vert \left\vert
\nabla \right\vert ^{\alpha }\phi (x,t)\right\vert ^{2}}{\left\vert
x\right\vert ^{\beta -2\alpha }}dxdt \\
&=&\frac{1}{\beta }\frac{1}{(2\pi )^{n}}\int_{\mathbb{S}^{n-1}}dS_{\omega
_{1}}\int_{\mathbb{S}^{n-1}}dS_{\omega _{2}}\int_{\mathbb{R}^{n}}\frac{%
e^{ix\cdot r_{1}(\omega _{1}-\omega _{2})}}{\left\vert x\right\vert ^{\beta
-2\alpha }}dx \\
&&\int_{0}^{\infty }\left( r_{1}^{n-\beta +2\alpha }\hat{\phi}%
_{0}(r_{1}\omega _{1})\overline{\hat{\phi}_{0}(r_{1}\omega _{2})}\right)
r_{1}^{n-1}dr_{1} \\
&=&c_{n,\alpha }\int_{0}^{\infty }r_{1}^{n-1}dr_{1}\int_{\mathbb{S}%
^{n-1}}dS_{\omega _{1}}\int_{\mathbb{S}^{n-1}}dS_{\omega _{2}}\frac{1}{%
\left\vert \omega _{1}-\omega _{2}\right\vert ^{n-\beta +2\alpha }}\hat{\phi}%
_{0}(r_{1}\omega _{1})\overline{\hat{\phi}_{0}(r_{1}\omega _{2})},
\end{eqnarray*}%
excluding the case when $\beta -2\alpha =n$ due to the fact that $\left\vert
x\right\vert ^{-n}$ is not a tempered distribution in $n$ d.

Because $n-\beta +2\alpha <n-1$ if $1<\beta -2\alpha ,$ the above
computation concludes the proof of theorem \ref%
{FullKatoEstimateForFreeSchrodinger}.

\begin{remark}
The steps in the above proof can be traced back to Sj\"{o}lin \cite{Sjolin}
in which the author proved various other local smoothing estimates for the
free Schr\"{o}dinger equation. In the case we are dealing with here, the
computation is carried out explicitly.
\end{remark}

\subsection{Proof of theorems \protect\ref{1d odd estimate for Schrodinger}
and \protect\ref{1/2derivative for 1d}}

Relation \ref{integrate in time} reads%
\begin{eqnarray*}
&&\int_{-\infty }^{\infty }\left\vert \left\vert \nabla \right\vert ^{\alpha
}\phi (x,t)\right\vert ^{2}dt \\
&=&\frac{1}{4\pi ^{2}}\int_{-\infty }^{\infty }dt\int_{-\infty }^{\infty
}d\xi _{1}\int_{-\infty }^{\infty }d\xi _{2}\left( e^{ix\xi
_{1}}e^{-i\left\vert \xi _{1}\right\vert ^{\beta }t}e^{-ix\xi
_{2}}e^{i\left\vert \xi _{2}\right\vert ^{\beta }t}\left\vert \xi
_{1}\right\vert ^{\alpha }\left\vert \xi _{2}\right\vert ^{\alpha }\hat{\phi}%
_{0}(\xi _{1})\overline{\hat{\phi}_{0}(\xi _{2})}\right) \\
&=&\frac{1}{4\pi ^{2}}\int_{-\infty }^{\infty }d\xi _{1}\lim_{\epsilon
\rightarrow 0}\int_{-\infty }^{\infty }d\xi _{2}\left( e^{ix(\xi _{1}-\xi
_{2})}\hat{\eta}_{\varepsilon }(\left\vert \xi _{1}\right\vert ^{\beta
}-\left\vert \xi _{2}\right\vert ^{\beta })\left\vert \xi _{1}\right\vert
^{\alpha }\left\vert \xi _{2}\right\vert ^{\alpha }\hat{\phi}_{0}(\xi _{1})%
\overline{\hat{\phi}_{0}(\xi _{2})}\right) \\
&=&\frac{1}{4\pi ^{2}}\int_{0}^{\infty }d\xi _{1}\lim_{\epsilon \rightarrow
0}\int_{0}^{\infty }d\xi _{2}+\frac{1}{4\pi ^{2}}\int_{-\infty }^{0}d\xi
_{1}\lim_{\epsilon \rightarrow 0}\int_{-\infty }^{0}d\xi _{2} \\
&&+\frac{1}{4\pi ^{2}}\int_{0}^{\infty }d\xi _{1}\lim_{\epsilon \rightarrow
0}\int_{-\infty }^{0}d\xi _{2}+\frac{1}{4\pi ^{2}}\int_{-\infty }^{0}d\xi
_{1}\lim_{\epsilon \rightarrow 0}\int_{0}^{\infty }d\xi _{2}
\end{eqnarray*}%
With the same procedure in the proof of theorem \ref%
{FullKatoEstimateForFreeSchrodinger}, we deduce%
\begin{eqnarray*}
&&\int_{-\infty }^{\infty }\frac{1}{\left\vert x\right\vert ^{2-2\alpha }}%
dx\int_{-\infty }^{\infty }\left\vert \left\vert \nabla \right\vert ^{\alpha
}\phi (x,t)\right\vert ^{2}dt \\
&=&\frac{1}{2\pi }\int_{-\infty }^{\infty }\frac{1}{\left\vert x\right\vert
^{\beta -2\alpha }}dx(\int_{0}^{\infty }\frac{1}{\beta \left\vert \xi
_{1}\right\vert ^{\beta -1-2\alpha }}\hat{\phi}_{0}(\xi _{1})\overline{\hat{%
\phi}_{0}(\xi _{1})}d\xi _{1}+\int_{-\infty }^{0}\frac{1}{\beta \left\vert
\xi _{1}\right\vert ^{\beta -1-2\alpha }}\hat{\phi}_{0}(\xi _{1})\overline{%
\hat{\phi}_{0}(\xi _{1})}d\xi _{1} \\
&&+\int_{0}^{\infty }\frac{e^{2ix\xi _{1}}}{\beta \left\vert \xi
_{1}\right\vert ^{\beta -1-2\alpha }}\hat{\phi}_{0}(\xi _{1})\overline{\hat{%
\phi}_{0}(-\xi _{1})}d\xi _{1}+\int_{-\infty }^{0}\frac{e^{2ix\xi _{1}}}{%
\beta \left\vert \xi _{1}\right\vert ^{\beta -1-2\alpha }}\hat{\phi}_{0}(\xi
_{1})\overline{\hat{\phi}_{0}(-\xi _{1})}d\xi _{1}) \\
&=&\frac{1}{2\pi }\int_{-\infty }^{\infty }\frac{1}{\left\vert x\right\vert
^{\beta -2\alpha }}\int_{-\infty }^{\infty }\frac{1-e^{2ix\xi _{1}}}{\beta
\left\vert \xi _{1}\right\vert ^{\beta -1-2\alpha }}\left\vert \hat{\phi}%
_{0}(\xi _{1})\right\vert ^{2}d\xi _{1} \\
&=&\frac{1}{2\beta \pi }\int_{-\infty }^{\infty }d\xi _{1}\frac{\left\vert 
\hat{\phi}_{0}(\xi _{1})\right\vert ^{2}}{\left\vert \xi _{1}\right\vert
^{\beta -1-2\alpha }}\int_{-\infty }^{\infty }dx\frac{1-\cos 2x\xi _{1}}{%
\left\vert x\right\vert ^{\beta -2\alpha }}
\end{eqnarray*}%
becasue $\hat{\phi}_{0}$ is odd if $\phi _{0}$ is odd. However,%
\begin{eqnarray*}
&&\int_{-\infty }^{\infty }\frac{1-\cos 2x\xi _{1}}{\left\vert x\right\vert
^{\beta -2\alpha }}dx \\
&=&2\int_{0}^{\infty }\frac{1-\cos 2x\xi _{1}}{x^{\beta -2\alpha }}dx \\
&=&\frac{2\cdot 2\xi _{1}}{\beta -1-2\alpha }\int_{0}^{\infty }\frac{\sin
2x\xi _{1}}{x^{\beta -1-2\alpha }}dx \\
&=&\frac{2\cdot 2\left\vert \xi _{1}\right\vert }{\beta -1-2\alpha }\frac{%
\Gamma (2-\beta +2\alpha )\sin \left( \frac{2-\beta +2\alpha }{2}\right) \pi 
}{(2\left\vert \xi _{1}\right\vert )^{2-\beta +2\alpha }} \\
&=&\frac{2^{\beta -2\alpha }\Gamma (2-\beta +2\alpha )\sin \left( \frac{%
2-\beta +2\alpha }{2}\pi \right) }{\beta -1-2\alpha }\left\vert \xi
_{1}\right\vert ^{\beta -1-2\alpha }
\end{eqnarray*}%
valid when $1<\beta -2\alpha \leqslant 2$ i.e.%
\begin{eqnarray*}
\int_{-\infty }^{\infty }\int_{-\infty }^{\infty }\frac{\left\vert
\left\vert \nabla \right\vert ^{\alpha }\phi (x,t)\right\vert ^{2}}{%
\left\vert x\right\vert ^{\beta -2\alpha }}dxdt &=&\frac{2^{\beta -2\alpha
}\Gamma (2-\beta +2\alpha )\sin \left( \frac{2-\beta +2\alpha }{2}\pi
\right) }{2\beta (\beta -1-2\alpha )\pi }\left\Vert \hat{\phi}%
_{0}\right\Vert _{2}^{2} \\
&=&\frac{2^{\beta -2\alpha }\Gamma (2-\beta +2\alpha )\sin \left( \frac{%
2-\beta +2\alpha }{2}\pi \right) }{\beta (\beta -1-2\alpha )}\left\Vert \phi
_{0}\right\Vert _{2}^{2}.
\end{eqnarray*}%
So theorem \ref{1d odd estimate for Schrodinger} is concluded. Notice that
relation \ref{integrate in time} becomes%
\begin{equation*}
\int_{-\infty }^{\infty }\left\vert \left\vert \nabla \right\vert ^{\alpha
}\phi (x,t)\right\vert ^{2}dt=\frac{1}{2\pi }\int_{-\infty }^{\infty }\frac{%
1+e^{2ix\xi _{1}}}{\beta \left\vert \xi _{1}\right\vert ^{\beta -1-2\alpha }}%
\left\vert \hat{\phi}_{0}(\xi _{1})\right\vert ^{2}d\xi _{1}
\end{equation*}%
if the initial data $\phi _{0}$ is even. Via the odd-even decomposition, we
have also proven theorem \ref{1/2derivative for 1d}.

\end{document}